\documentclass[preprint]{elsarticle}
\usepackage{amsfonts}
\usepackage{amsmath}
\usepackage[all]{xy}
\usepackage{tikz}
\usetikzlibrary{trees}
\usetikzlibrary{arrows,shapes,snakes,automata,backgrounds,petri}
\newtheorem{thm}{Theorem}
\newtheorem{lem}[thm]{Lemma}
\newtheorem{cor}[thm]{Corollary}
\newproof{pf}{Proof}
\newproof{pot1}{Proof of Theorem \ref{main}}
\newproof{pot2}{Proof of Lemma \ref{lemconv}}
\newproof{pot3}{Proof of Corollary \ref{corconv}}
\newproof{pot4}{Proof of Theorem \ref{corconv1}}

\begin{document}
\begin{frontmatter}

\title{$P$-matrices and signed digraphs}

\author[ref1]{Murad Banaji\corref{cor1}}
\author[ref2]{Carrie Rutherford}

\address[ref1]{Department of Mathematics, University College London, Gower Street, London WC1E 6BT, UK.}
\address[ref2]{Department of Business Studies, London South Bank University, 103 Borough Road, London SE1 0AA, UK}

\cortext[cor1]{Corresponding author: m.banaji@ucl.ac.uk. }

\begin{abstract}
We associate a signed digraph with a list of matrices whose dimensions permit them to be multiplied, and whose product is square. Cycles in this graph have a parity, that is, they are either even (termed e-cycles) or odd (termed o-cycles). The absence of e-cycles in the graph is shown to imply that the matrix product is a $P_0$-matrix, i.e., all of its principal minors are nonnegative. Conversely, the presence of an e-cycle is shown to imply that there exists a list of matrices associated with the graph whose product fails to be a $P_0$-matrix. The results generalise a number of previous results relating $P$- and $P_0$-matrices to graphs. 
\end{abstract}


\end{frontmatter}

\section{Introduction and statement of the main results}

$P$-matrices are square matrices, all of whose principal minors are positive. $P_0$-matrices \cite{hershkowitz} are square matrices all of whose principal minors are nonnegative, i.e., matrices in the closure of the $P$-matrices. We will be interested in real $P_0$-matrices. A well-known result of Gale and Nikaido \cite{gale}, whose extensions and corollaries are discussed in \cite{parthasarathy}, states that a differentiable function on a rectangular domain in $\mathbb{R}^n$ with $P$-matrix Jacobian is injective. This result has a number of practical applications -- see for example \cite{banajiSIAM,donnell}. 

The results of Gale and Nikaido have various graph-theoretic implications, explored in \cite{soule,banajicraciun2,banajiJMAA} for example. The key results in these papers involve determining sufficient graph-theoretic conditions for a set of matrices to be $P_0$-matrices, followed by additional ``nondegeneracy'' conditions which guarantee that they are in fact $P$-matrices. Here, a general result is developed, of which some of these previous results become corollaries.

Let $k$ be an arbitrary positive integer. Throughout this paper, a subscript or superscript $j$ assumed to belong to $\{0, \ldots, k-1\}$ should be read as $j \bmod k$. Let $n_0, \ldots, n_{k-1}$ be positive integers. For each $j \in \{0, \ldots, k-1\}$, let $A^{(j)}$ be an $n_j \times n_{j+1}$ matrix, and define the $n_0 \times n_0$ matrix $A= A^{(0)}A^{(1)}\cdots A^{(k-1)}$. We will associate with $[A^{(0)}, A^{(1)},\ldots, A^{(k-1)}]$ a signed digraph $G_{A^{(0)}\cdots A^{(k-1)}}$, which will belong to a category of graphs termed signed $(k,\{1\})$-BC digraphs, to be defined below. The structure of these graphs will imply that all its cycles have length which is a multiple of $k$.  

Given a cycle $C$ with $kr_1$ edges, of which $r_2$ have negative sign, we define $C$ to be an {\bf e-cycle} if $(-1)^{r_1+r_2}=1$ and an {\bf o-cycle} otherwise. A signed $(k,\{1\})$-BC digraph containing no e-cycles will be termed ``e-cycle-free''. The first main theorem in this paper is:

\begin{thm}
\label{main}
If $G = G_{A^{(0)}\cdots A^{(k-1)}}$ is e-cycle-free then $A = A^{(0)}A^{(1)}\cdots A^{(k-1)}$ is a $P_0$-matrix.
\end{thm}

A matrix $M$ determines the qualitative class $\mathcal{Q}(M)$ \cite{brualdi} consisting of all matrices with the same sign pattern as $M$. Explicitly, $\mathcal{Q}(M)$ consists of all matrices $X$ with the same dimensions as $M$, and satisfying $M_{ij} > 0 \Rightarrow X_{ij} > 0$, $M_{ij} < 0 \Rightarrow X_{ij} < 0$ and $M_{ij} = 0 \Rightarrow X_{ij} = 0$. Given two matrices $M$ and $N$ of dimensions such that they can be multiplied, we write 
\[
\mathcal{Q}(M)\mathcal{Q}(N) = \{M^{'}N^{'}\,|\,M^{'} \in \mathcal{Q}(M), N^{'} \in \mathcal{Q}(N)\}. 
\]
This definition extends naturally to any ordered set of multiplicable matrices. Note that $\mathcal{Q}(M)\mathcal{Q}(N)$ is in general different from $\mathcal{Q}(MN)$, and is often not a subset of any qualitative class. The second main theorem in this paper is:

\begin{thm}
\label{corconv1}
All matrices in $\mathcal{Q}(A^{(0)})\mathcal{Q}(A^{(1)})\cdots \mathcal{Q}(A^{(k-1)})$ are $P_0$-matrices if and only if $G_{A^{(0)}A^{(1)}\cdots A^{(k-1)}}$ is e-cycle-free. 
\end{thm}

\section{Signed $(k,\{1\})$-BC digraphs}

Consider a digraph $G$ with vertex set $V(G)$ and edge set $E(G)$. Let $S$ be any set of integers. $G$ will be termed ``$(k,S)$-block circulant'', abbreviated to $(k,S)$-BC, if
\begin{enumerate}
\item $V(G)$ is partitioned into $k$ sets $\{V_0, V_1, \ldots, V_{k-1}\}$. 
\item For $k \geq 2$, if $(j - i) \bmod k \not \in S$, then there is no (directed) edge from a vertex in $V_i$ to a vertex in $V_j$. 
\end{enumerate}
Every digraph is trivially a $(1,S)$-BC digraph for arbitrary $S$. When $k \geq 2$, $(k,S)$-BC digraphs are a generalisation of circulant digraphs (\cite{alspach,boesch} for example). Note, however, that vertices in a $(k,S)$-BC digraph may have arbitrary outdegree and indegree. Here, only the special case $S = \{1\}$ concerns us. In a $(k,\{1\})$-BC digraph, a (directed) path from a vertex in $V_j$ to a vertex in $V_j$ must include a vertex from each $V_{j^{'}}$, $j^{'} \not = j$. It follows that all cycles in a $(k,\{1\})$-BC digraph have length which is a multiple of $k$.

{\bf Remark.} Although, for $k \geq 2$, $(k,\{1\})$-BC digraphs are $k$-colourable, $k$ is not in general the chromatic number of $G$: for example, a $(2r,\{1\})$-BC digraph with nonempty edge-set is in fact bipartite.\\

Let $V_j$ contain $n_j$ vertices. Assume some ordering on these vertices and let $V_j^i$ ($i \in \{1, \ldots, n_j\}$) refer to the $i$th vertex in $V_j$. As usual, an edge $(v, \tilde v)$ refers to the edge directed from $v$ to $\tilde v$. 

A digraph $G$ is signed if there is a function $\mathrm{sign}:E(G) \to \{-1, +1\}$. A signed $(k,\{1\})$-BC digraph $G = G_{A^{(0)}\cdots A^{(k-1)}}$ is associated with a list of matrices $[A^{(0)}, \ldots, A^{(k-1)}]$ as above in a very simple way: (i) for each $j = 0, \ldots, k-1$, $|V_j| = n_j$, (ii) there exists an edge $(V_j^{r},V_{j + 1}^{s})$ in $G$ if and only if $(A^{(j)})_{rs} \not = 0$, and (iii) the edge $(V_j^{r},V_{j + 1}^{s})$ takes the sign of $(A^{(j)})_{rs}$. Note that entries in the matrices $A^{(j)}$ are in one-to-one correspondence with edges in $G$, and that the sign-pattern of $A^{(j)}$ is in fact a block in the (signed) adjacency matrix of $G$.\\

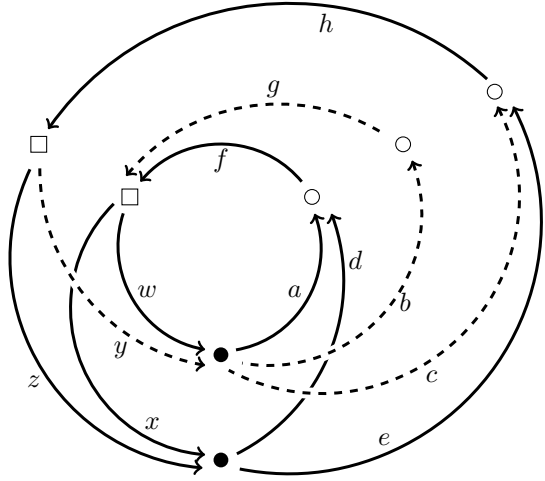
\begin{figure}[h]

\begin{center}
\begin{tikzpicture}[domain=0:4,scale=0.7]

\fill (0,-2) circle (4pt);
\fill (0,-4) circle (4pt);

\draw (-1.88,0.85) rectangle (-1.58,1.15);
\draw (-3.61,1.85) rectangle (-3.31,2.15);

\draw (1.73,1) circle (4pt);
\draw (3.46,2) circle (4pt);
\draw (5.20,3) circle (4pt);

\draw[<-, line width=0.04cm] (-0.3,-1.9) arc (260:162:2cm);
\draw[<-, line width=0.04cm] (-0.35,-4.15) arc (265:155:4cm);

\draw[<-, line width=0.04cm] (-0.3,-3.9) arc (265:135:2.8cm);
\draw[-, line width=0.15cm, color=white] (-0.3,-2.2) arc (260:186:3.8cm);
\draw[<-, dashed, line width=0.04cm] (-0.3,-2.2) arc (260:180:3.8cm);

\draw[<-, line width=0.04cm] (-1.53,1.29) arc (140:40:2cm);
\draw[<-, line width=0.04cm] (-3.28,2.29) arc (145:48:5.6cm);
\draw[<-, dashed, line width=0.04cm] (-1.8,1.4) arc (140:60:3.8cm);

\draw[<-, line width=0.04cm] (1.78, 0.68) arc (20:-80:1.95cm);
\draw[<-, line width=0.04cm] (2.1, 0.68) arc (20:-63:3.7cm);
\draw[-, line width=0.15cm,color=white] (3.63, 1.69) arc (20:-90:2.9cm);
\draw[<-, dashed, line width=0.04cm] (3.63, 1.69) arc (20:-100:2.9cm);

\draw[-, line width=0.15cm, color=white] (5.25, 2.72) arc (27:-110:3.75cm);
\draw[<-, dashed, line width=0.04cm] (5.25, 2.72) arc (27:-118:3.75cm);

\draw[<-, line width=0.04cm] (5.55, 2.72) arc (27:-101:4.8cm);

\node at (-1.4,-0.8) {$w$};
\node at (-1.3,-3.3) {$x$};
\node at (-1.9,-1.9) {$y$};
\node at (-3.55,-2.5) {$z$};

\node at (0,1.7) {$f$};
\node at (1,3.05) {$g$};
\node at (2,4.3) {$h$};

\node at (1.4,-0.8) {$a$};
\node at (2.55,-0.2) {$d$};
\node at (3.5,-1) {$b$};
\node at (4,-2.4) {$c$};
\node at (3.1,-3.6) {$e$};

\end{tikzpicture}
\end{center}
\caption{\label{3SR}The signed $(3,\{1\})$-BC digraph corresponding to the product of the three matrices $A^{(0)}A^{(1)}A^{(2)}$ in Eq.~\ref{examprod}. The graph has been laid out to emphasise its circulant structure. The vertices in $V_0, V_1$ and $V_2$ are represented as filled circles, open circles and boxes respectively. Negative edges are represented as dashed lines while positive edges are bold lines. Labels on each edge represent the absolute values of the corresponding entries in the matrices, and are not strictly part of the graph: they have been added to indicate the correspondence between edges and matrix entries.}
\end{figure}

{\bf Example.} As an example consider the matrices:
\begin{equation}
\label{examprod}
A^{(0)} = \left(\begin{array}{rrr}a&-b&-c\\d&0&e\end{array}\right)\,,\quad A^{(1)} = \left(\begin{array}{rr}f&0\\-g&0\\0&h\end{array}\right)\,,\quad A^{(2)} = \left(\begin{array}{rr}w&x\\-y&z\end{array}\right)\,,
\end{equation}
where $a,b,c,d,e,f,g,h,w,x,y$ and $z$ are arbitrary positive real numbers. Associated with the product $A^{(0)}A^{(1)}A^{(2)}$ is the signed $(3,\{1\})$-BC digraph shown in Figure~\ref{3SR}. 

Defining $\mathbf{0}_{m \times n}$ to be the $m \times n$ zero matrix, the signed adjacency matrix of the graph $G$ in Figure~\ref{3SR} has block structure:
\[
\left(\begin{array}{ccc}
\mathbf{0}_{2 \times 2} & \left(\begin{array}{rrr}1&-1&-1\\1&0&1\end{array}\right) & \mathbf{0}_{2 \times 2} \\
 \mathbf{0}_{3 \times 2} & \mathbf{0}_{3 \times 3} & \left(\begin{array}{rr}1&0\\-1&0\\0 & 1\end{array}\right)\\
\left(\begin{array}{rr}1&1\\-1&1\end{array}\right) & \mathbf{0}_{2 \times 3}& \mathbf{0}_{2 \times 2}
\end{array}\right)\,.
\]
It can be seen immediately that each block is simply the sign-pattern of $A^{(0)}$, $A^{(1)}$ or $A^{(2)}$. 

Although the graph in Figure~\ref{3SR} has a number of cycles (both of length $3$ and of length $6$), all of these can be computed to be o-cycles, and so, by Theorem~\ref{main}, the product $A^{(0)}A^{(1)}A^{(2)}$ is a $P_0$-matrix. This is true whatever the magnitudes of the entries in the matrices. Clearly, given isomorphic signed digraphs $G_1 \cong G_2$, $G_1$ is e-cycle-free if and only if $G_2$ is e-cycle-free. Since 
\[
G_{A^{(1)}A^{(2)}A^{(0)}} \cong G_{A^{(2)}A^{(0)}A^{(1)}} \cong G_{A^{(0)}A^{(1)}A^{(2)}},
\] 
$A^{(1)}A^{(2)}A^{(0)}$ and $A^{(2)}A^{(0)}A^{(1)}$ are also $P_0$-matrices. 

\section{Preliminaries needed for the proofs}

{\bf Permutations.} Consider an ordered set $\alpha = [\alpha_1, \alpha_2, \ldots, \alpha_r]$, and a bijection $\beta:\alpha \to \alpha$. Defining $\beta_i \equiv \beta(\alpha_i)$, the ordered set $[\beta_1, \beta_2, \ldots, \beta_r]$ will also be referred to as $\beta$. It will always be clear from context whether an object referred to is a bijection or an ordered set.

Any permutation $\beta$ has a parity $P(\beta)$, i.e., $P(\beta) = +1$ if $\beta$ is an even permutation and $P(\beta) = -1$ otherwise. Given two permutations $\beta$ and $\beta^{'}$, $P(\beta\beta^{'}) = P(\beta)P(\beta^{'})$, implying that $P(\beta) = P(\beta^{-1})$. Note the following elementary result about the parity of permutations.
\begin{lem}
\label{permsigns}
Consider a permutation $\beta$ of a finite set of size $r$. Write $\beta$ as a product of disjoint cycles, $C_1, \ldots, C_s$ ($1 \leq s \leq r$), including trivial cycles. Then 
\[
P(\beta) = (-1)^{r - s}
\]
i.e., $\beta$ is even (resp. odd) if the total number of elements in $\beta$ minus the total number of cycles in its decomposition is even (resp. odd). 
\end{lem}
\begin{pf}
See \cite{banajicraciun}, for example.
\end{pf}

From here on $\alpha^{(j)}$ will always refer to a nonempty subset of $\{1, \ldots, n_j\}$, and will be assumed to have the natural ordering. $\alpha^{(j)}_{m}$ will refer to the $m$th element in $\alpha^{(j)}$ so that $\alpha^{(j)}_{1} < \alpha^{(j)}_{2} < \alpha^{(j)}_{3} < \cdots$. Given some $\alpha^{(j)}$, define $V_j^{\alpha^{(j)}} = \{V_j^k\,|\,k \in \alpha^{(j)}\} \subseteq V_j$. $\beta^{(j)}$ will refer to a permutation of $\alpha^{(j)}$. Given the one-to-one correspondence between the elements in $\alpha^{(j)}$, and vertices in $V_j^{\alpha^{(j)}}$, $\beta^{(j)}$ can equally be regarded as a permutation on $V_j^{\alpha^{(j)}}$.

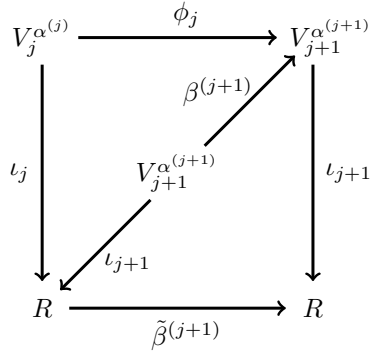
\begin{figure}[h]

\begin{center}
\begin{tikzpicture}[domain=0:4,scale=1.2]


\node at (1,4) {$V_{j}^{\alpha^{(j)}}$};
\node at (4.2,4) {$V_{j+1}^{\alpha^{(j+1)}}$};
\node at (1,1) {$R$};
\node at (4,1) {$R$};
\node at (2.5,2.5) {$V_{j+1}^{\alpha^{(j+1)}}$};

\draw[->, line width=0.04cm] (1.4,4.0) -- (3.6, 4.0);
\draw[->, line width=0.04cm] (1.3,1.0) -- (3.7, 1.0);

\draw[->, line width=0.04cm] (1.0, 3.7) -- (1.0, 1.3);
\draw[->, line width=0.04cm] (4.0, 3.7) -- (4.0, 1.3);

\draw[<-, line width=0.04cm] (1.2, 1.2) -- (2.2, 2.2);
\draw[->, line width=0.04cm] (2.8, 2.8) -- (3.8, 3.8);

\node at (0.75, 2.5) {$\iota_{j}$};
\node at (2.95, 3.4) {$\beta^{(j+1)}$};
\node at (2.6, 4.25) {$\phi_j$};
\node at (4.4, 2.5) {$\iota_{j+1}$};

\node at (2.6, 0.7) {$\tilde \beta^{(j+1)}$};
\node at (1.95, 1.5) {$\iota_{j+1}$};
\end{tikzpicture}
\end{center}
\caption{\label{CD2}The commutative diagram which encapsulates the relationships $\tilde \beta^{(j)} = \iota_j\circ\beta^{(j)}\circ \iota_j^{-1}$, $\phi_j = \beta^{(j+1)}\circ \iota_{j+1}^{-1}\circ \iota_{j}$, and $\phi_j = \iota_{j+1}^{-1}\circ \tilde \beta^{(j+1)}\circ \iota_{j}$.}
\end{figure}

Now consider some sequence $(\alpha^{(0)}, \ldots, \alpha^{(k-1)})$, such that $|\alpha^{(i)}| = |\alpha^{(j)}|$ for each $i,j$, and a corresponding sequence of permutations $(\beta^{(0)}, \ldots, \beta^{(k-1)})$. Define $R=[1, 2, \ldots, |\alpha^{(0)}|]$, and define the bijections $\iota_j:V_j^{\alpha^{(j)}} \to R$ by $\iota_j(V_j^{\alpha^{(j)}_{m}}) = m$. In other words, $\iota_j$ associates with each vertex in $V_j^{\alpha^{(j)}}$ its order. $\beta^{(j)}$ then induces the bijection $\tilde \beta^{(j)}:R \to R$ defined by $\tilde \beta^{(j)} = \iota_j\circ\beta^{(j)}\circ \iota_j^{-1}$. Further, define the bijections $\phi_j: V_{j}^{\alpha^{(j)}} \to V_{j+1}^{\alpha^{(j+1)}}$ by $\phi_j = \iota_{j+1}^{-1}\circ \tilde \beta^{(j+1)}\circ \iota_{j}$, i.e. $\phi_j(V_{j}^{\alpha^{(j)}_{m}}) = V_{j+1}^{\beta^{(j+1)}_{m}}$. Equivalently, $\phi_j = \beta^{(j+1)}\circ \iota_{j+1}^{-1}\circ \iota_{j}$. Figure~\ref{CD2} illustrates all of these relationships. Note that in the special case $k = 1$, $\phi_j = \beta^{(j)}$. 

\section{Proof of Theorem~\ref{main} and its immediate consequences}

The following notation is used. Given an $r \times s$ matrix $M$, and two (nonempty) ordered sets $\gamma \subseteq \{1, \ldots, r\}$ and $\delta \subseteq \{1,\ldots, s\}$, then $M(\gamma|\delta)$ is the submatrix of $M$ with rows indexed by $\gamma$ and columns indexed by $\delta$. If $|\gamma| = |\delta|$, then $M[\gamma|\delta] \equiv \mathrm{det}(M(\gamma|\delta))$. We write $M[\gamma]$ as shorthand for $M[\gamma|\gamma]$. If $\gamma$ and $\delta$ are of equal size, then $M_{\gamma,\delta}$ will refer to $\prod_{i=1}^{|\gamma|}M_{\gamma_i, \delta_i}$.\\

\begin{pot1}
The cases $k=1$ and $k \geq 2$ are conceptually similar: however in order to avoid notational difficulties, they are presented separately.\\

{\bf Case 1: $k=1$.} Let $V = V(G)$. Choose and fix some nonempty $\alpha \subseteq \{1, \ldots, n\}$, and consider the minor $A[\alpha]$. Given a permutation $\beta$ of $\alpha$, define $T$, a term in $A[\alpha]$, by:
\[
T = P(\beta)A_{\alpha, \beta}\,.
\]
Suppose that $T \not = 0$. Since nonzero entries in $A$ are in one-to-one correspondence with edges in $G$, corresponding to $T$ is an edge-set $E \subseteq E(G)$ defined as follows: the edge $(V^r, V^s)$ is in $E$ if and only if $(r,s) = (\alpha_{m}, \beta_{m})$ for some integer $m$. Equivalently, the edge $(V^r, V^s)$ is in $E$ if and only if $V^s = \beta(V^r)$.\\

The set of endpoints of edges in $E$ is precisely $V^{\alpha}$. Consider the vertex $V^{\alpha_{m}}$. Then the incoming edge  $\left(\beta^{-1}(V^{\alpha_{m}}), V^{\alpha_{m}}\right)$ and the outgoing edge $\left(V^{\alpha_{m}}, \beta(V^{\alpha_{m}})\right)$ are the only edges from $E$ incident on $V^{\alpha_{m}}$. If $\beta(V^{\alpha_{m}}) = V^{\alpha_{m}}$, then these edges coincide and in fact there is a loop at $V^{\alpha_{m}}$. Otherwise, the edges are distinct. In either case, $E$, regarded as a subgraph of $G$, consists of vertex- and edge-disjoint cycles.

Suppose $\beta^{m}(v) = v$ for some vertex $v$, but $\beta^{q}(v) \not = v$ for $q < m$. Then the vertex $v$ lies on a cycle in $E$ of length $m$. Thus, decomposing $\beta$ as a product of disjoint cycles (including trivial cycles), these cycles are in one-to-one correspondence with cycles -- in the graph-theoretic sense -- in $E$. Trivial cycles correspond to loops. Assume that there are $N$ such cycles. By Lemma~\ref{permsigns}, $P(\beta) = (-1)^{|\alpha|-N}$.

Returning to the term $T$ in the expansion of $A[\alpha]$, 
\begin{equation}
\label{eqpt0}
\mathrm{sign}(T) = P(\beta)\mathrm{sign}(A_{\alpha, \beta})
\end{equation}
Consider a cycle $C$ in $E$ including $r_1$ edges of which $r_2$ are negative, and define $\mathrm{sign}(C) = (-1)^{r_2}$ (i.e. $\mathrm{sign}(c)$ is the product of signs of edges in $C$). If $C$ is an e-cycle, then $(-1)^{r_1+r_2}=1$, and so $\mathrm{sign}(C) = (-1)^{r_1}$. Similarly if $C$ is an o-cycle, then $\mathrm{sign}(C) = (-1)^{r_1+1}$. Decompose $E$ into disjoint cycles, which comprise $N_e$ e-cycles, and $N_o$ o-cycles (so that $N = N_o + N_e$). Let $\theta$ be the total number of edges in o-cycles and $\theta_e$ the total number of edges in e-cycles, so that $\theta_o + \theta_e = |\alpha|$ (since there are $|\alpha|$ edges in $E$). Taking the product of signs of edges in $E$ over e-cycles and o-cycles separately gives
\begin{equation}
\label{eqpt2}
\mathrm{sign}\left(A_{\alpha, \beta}\right) = (-1)^{\theta_o}(-1)^{\theta_e + N_o} = (-1)^{|\alpha| + N_o}\,.
\end{equation}
Since $P(\beta) = (-1)^{|\alpha|-N}$, and $\mathrm{sign}\left(A_{\alpha, \beta}\right) = (-1)^{|\alpha| + N_o}$, Eq.~\ref{eqpt0} gives:
\begin{equation}
\label{eqpt3}
\mathrm{sign}(T) = (-1)^{|\alpha|-N}(-1)^{|\alpha| + N_o} = (-1)^{N_e}.
\end{equation}
If $G$ is e-cycle-free, then $N_e = 0$ in this expression, in which case $\mathrm{sign}(T) = 1$. Since $T$ is an arbitrary nonzero term in $A[\alpha]$, $A[\alpha] \geq 0$. Since $\alpha$ is an arbitrary nonempty subset of $\{1, \ldots, n\}$, $A$ is a $P_0$-matrix. \\

{\bf Case 2: $k \geq 2$.} Choose and fix some nonempty $\alpha^{(0)} \subseteq \{1, \ldots, n_0\}$, and consider the minor $A[\alpha^{(0)}]$. Enumerate all lists $(\alpha^{(1)}, \ldots, \alpha^{(k-1)})$ where for each $i = 1, \ldots, k-1$, $\alpha^{(i)}$ is a nonempty subset of $\{1, \ldots, n_i\}$ such that $|\alpha^{(i)}| = |\alpha^{(0)}|$. There may of course be no such subset for some $i$, and hence no such lists. 

By the Cauchy-Binet formula \cite{gantmacher} applied recursively, 

\[
A[\alpha^{(0)}] = \sum_{\substack{(\alpha^{(1)}, \ldots, \alpha^{(k-1)}),\\ |\alpha^{(i)}| = |\alpha^{(0)}|}}\left(\prod_{j=0}^{k-1} A^{(j)}[\alpha^{(j)}|\alpha^{(j+1)}]\right)\,.
\]
The sum is over all possible lists $(\alpha^{(1)}, \ldots, \alpha^{(k-1)})$ with $|\alpha^{(i)}| = |\alpha^{(0)}|$. Now choose and fix some particular choice $\alpha^{(1)}, \ldots, \alpha^{(k-1)}$, and choose permutations $\beta^{(0)}, \ldots, \beta^{(k-1)}$. For each $j$ define $T_j$, a term in $A^{(j)}[\alpha^{(j)}|\alpha^{(j+1)}]$, by:
\[
T_j = P(\beta^{(j+1)})A^{(j)}_{\alpha^{(j)}, \beta^{(j+1)}}\,.
\]

Suppose that for each $j$, $T_j$ is nonzero so that $T = \prod_jT_j \not = 0$. Note that $T$ is then a nonzero term in the expansion of $A[\alpha^{(0)}]$. Since nonzero entries in the matrices $A^{(j)}$ are in one-to-one correspondence with edges in $G$, corresponding to $T$ is an edge-set $E \subseteq E(G)$ defined as follows: the edge $(V_{j}^r, V_{j+1}^s)$ is in $E$ if and only if $(r,s) = (\alpha^{(j)}_{m}, \beta^{(j+1)}_{m})$ for some integer $m$. Equivalently, the edge $(V_{j}^r, V_{j+1}^s)$ is in $E$ if and only if $\phi_{j}(V_{j}^r) = V_{j+1}^s$.\\

The set of endpoints of edges in $E$ is precisely $\bigcup_j V_j^{\alpha^{(j)}}$, and in fact each such vertex has exactly two edges from $E$ incident on it, one incoming and one outgoing. For example, consider the vertex $V_j^{r}$, where $r \in \alpha^{(j)}$. Then the incoming edge  $\left(\phi_{j-1}^{-1}(V_j^{r}), V_j^{r}\right)$ and the outgoing edge $\left(V_j^{r}, \phi_{j}(V^{r})\right)$ are distinct edges in $E$, and are, by the definition of $E$, the only two edges in $E$ incident on $V_j^{r}$. As a consequence, $E$, regarded as a subgraph of $G$, consists of vertex- and edge-disjoint cycles.

Next, consider the bijection $\phi: V_{0}^{\alpha^{(0)}} \to V_{0}^{\alpha^{(0)}}$ defined by $\phi = \phi_{k-1}\circ \phi_{k-2}\circ \cdots \circ \phi_{0}$. Suppose $\phi^m(v) = v$ for some vertex $v$, but $\phi^{q}(v) \not = v$ for $q < m$. Then the vertex $v$ lies on a cycle in $E$ of length $km$. Decomposing $\phi$ as a product of disjoint cycles (including trivial cycles), these cycles are in one-to-one correspondence with cycles -- in the graph-theoretic sense -- in $E$. Assume that there are $N$ such cycles. By Lemma~\ref{permsigns}, $P(\phi) = (-1)^{|\alpha^{(0)}|-N}$. Applying the relations $\phi_j = \iota_{j+1}^{-1}\circ \tilde \beta^{(j+1)}\circ \iota_{j}$ gives
\[
\phi = \iota_{0}^{-1}\circ \tilde \beta^{(0)}\circ \tilde \beta^{(k-1)}\circ\cdots \circ \tilde \beta^{(2)}\circ \tilde \beta^{(1)}\circ\iota_{0}
\]
so that $P(\phi) = \prod_{j=0}^{k-1}P(\tilde \beta^{(j)}) = \prod_{j=0}^{k-1}P(\beta^{(j)})$. \\

Returning to the term $T$ in the expansion of $A[\alpha^{(0)}]$, 
\begin{equation}
\label{eqpt01}
\mathrm{sign}(T) = \prod_{j=0}^{k-1}\mathrm{sign}(T_j) = \left(\prod_{j=0}^{k-1}P(\beta^{(j+1)})\right)\left(\prod_{j=0}^{k-1}\mathrm{sign}\left(A^{(j)}_{\alpha^{(j)}, \beta^{(j+1)}}\right)\right)\,.
\end{equation}

The first term in this expression has already been determined: from above,  
\begin{equation}
\label{eqpt11}
\prod_{j=0}^{k-1}P(\beta^{(j+1)}) = \prod_{j=0}^{k-1}P(\beta^{(j)}) = P(\phi) = (-1)^{|\alpha^{(0)}|-N}\,.
\end{equation}
Consider a cycle $C$ in $E$ including $kr_1$ edges of which $r_2$ are negative. As in the case $k=1$, if $C$ is an e-cycle, then $\mathrm{sign}(C) = (-1)^{r_1}$, while if $C$ is an o-cycle, then $\mathrm{sign}(C) = (-1)^{r_1+1}$. Decompose $E$ into disjoint cycles, which comprise $N_e$ e-cycles, and $N_o$ o-cycles. Let $k\theta_0$ be the total number of edges in o-cycles and $k\theta_e$ the total number of edges in e-cycles, so that $\theta_o + \theta_e = |\alpha^{(0)}|$ (since there are $k|\alpha^{(0)}|$ edges in $E$). Taking the product of signs of edges in $E$ over e-cycles and o-cycles separately gives
\begin{equation}
\label{eqpt21}
\prod_{j=0}^{k-1}\mathrm{sign}\left(A^{(j)}_{\alpha^{(j)}, \beta^{(j+1)}}\right) = (-1)^{\theta_o}(-1)^{\theta_e + N_o} = (-1)^{|\alpha^{(0)}| + N_o}\,.
\end{equation}
Substituting Eqs.~\ref{eqpt11}~and~\ref{eqpt21} into Eq.~\ref{eqpt01} gives:
\begin{equation}
\label{eqpt31}
\mathrm{sign}(T) = (-1)^{|\alpha^{(0)}|-N}(-1)^{|\alpha^{(0)}| + N_o} = (-1)^{N_e}.
\end{equation}
Note that this is just Eq.~\ref{eqpt3} again. As in the case $k=1$, if $G$ is e-cycle-free, then $\mathrm{sign}(T) = 1$, and since $T$ is an arbitrary nonzero term in $A[\alpha^{(0)}]$, $A[\alpha^{(0)}] \geq 0$. Since $\alpha^{(0)}$ is an arbitrary nonempty subset of $\{1, \ldots, n_0\}$, $A$ is a $P_0$-matrix. \qquad \qed
\end{pot1}

We have the following corollary to Theorem~\ref{main}:

\begin{cor}
\label{cormain}
Consider a square matrix $A= A^{(0)}A^{(1)}\cdots A^{(k-1)}$ such that $G = G_{A^{(0)}A^{(1)}\cdots A^{(k-1)}}$ is e-cycle-free. Then matrices in $\mathcal{Q}(A^{(0)})\mathcal{Q}(A^{(1)})\cdots \mathcal{Q}(A^{(k-1)})$ are all $P_0$-matrices. 
\end{cor}
\begin{pf}
By definition, any matrix $B \in \mathcal{Q}(A^{(0)})\mathcal{Q}(A^{(1)})\cdots \mathcal{Q}(A^{(k-1)})$ can be written $B= B^{(0)}B^{(1)}\cdots B^{(k-1)}$, where $B^{(j)} \in \mathcal{Q}(A^{(j)})$. But 
\[
G_{B^{(0)}B^{(1)}\cdots B^{(k-1)}} \cong G_{A^{(0)}A^{(1)}\cdots A^{(k-1)}},
\]
so, by Theorem~\ref{main}, $B$ is a $P_0$-matrix. \qquad \qed
\end{pf}

{\bf Remark.} For each $r = 1, \ldots, k-1$ the matrix product 
\[
A^{(r)}A^{(r+1)}\cdots A^{(k-1)}A^{(0)}\cdots A^{(r-1)}
\] 
gives rise to a graph isomorphic to $G=G_{A^{(0)}A^{(1)}\cdots A^{(k-1)}}$. Thus, in fact, if $G$ is e-cycle-free, then all matrices in 
\[
\mathcal{Q}(A^{(r)})\mathcal{Q}(A^{(r+1)})\cdots \mathcal{Q}(A^{(k-1)})\mathcal{Q}(A^{(0)})\cdots \mathcal{Q}(A^{(r-1)})
\] 
are $P_0$-matrices.

\section{Converse results and proof of Theorem~\ref{corconv1}}

A variety of converse results are possible, that is, results which guarantee that if a signed $(k,\{1\})$-BC graph contains e-cycles, then there exist matrices in some set which fail to be $P_0$-matrices. The most useful formulations depend on the application. Lemma~\ref{lemconv} is the basic result from which such results follow:

\begin{lem}
\label{lemconv}
Consider a graph $G = G_{A^{(0)}A^{(1)}\cdots A^{(k-1)}}$ such that all edges in $G$ lie on a single e-cycle $C$ of length $kr$. Then $A = A^{(0)}A^{(1)}\cdots A^{(k-1)}$ is not a $P_0$-matrix. 
\end{lem}
\begin{pf}
Define the sets $(\alpha^{(0)}, \ldots, \alpha^{(k-1)})$ by the stipulation that $s \in \alpha^{(i)}$ if and only if $V_i^s$ lies on $C$. Clearly $|\alpha^{(i)}| = r$ for each $i$. Each vertex in $V_j^{\alpha^{(j)}}$ lies on $C$ and hence has exactly two edges, one incoming, and one outgoing, incident on it. So it is possible to define bijections $\phi_j: V_{j}^{\alpha^{(j)}} \to V_{j+1}^{\alpha^{(j+1)}}$ as follows: given vertices $v \in V_{j}^{\alpha^{(j)}}$ and $\tilde v \in V_{j+1}^{\alpha^{(j+1)}}$, $\phi_j(v) = \tilde v$ if there is a directed edge $(v, \tilde v)$ in $E$. Each $\phi_j$ induces a permutation $\beta^{(j+1)}:V_{j+1}^{\alpha^{(j+1)}} \to V_{j+1}^{\alpha^{(j+1)}}$ defined by $\beta^{(j+1)} = \phi_{j}\circ \iota_{j}^{-1}\circ \iota_{j+1}$ (see Figure~\ref{CD2}). 

Consider the minor $A[\alpha^{(0)}]$. There is a nonzero term in this minor
\[
T = \prod_{j=0}^{k-1}T_j = \prod_{j=0}^{k-1}P(\beta^{(j+1)})A^{(j)}_{\alpha^{(j)}, \beta^{(j+1)}}\,.
\]
Moreover $T$ is the unique nonzero term in $A[\alpha^{(0)}]$: any other nonzero term would imply the existence of an index $j$ and a permutation $\delta:V_{j+1}^{\alpha^{(j+1)}} \to V_{j+1}^{\alpha^{(j+1)}}$, $\delta \not = \beta^{(j+1)}$, such that $A^{(j)}_{\alpha^{(j)}, \delta} \not = 0$. Letting $s$ be an index such that $\delta_{s} \not = \beta^{(j+1)}_{s}$, $A^{(j)}_{\alpha^{(j)}_{s}, \delta_{s}}$ must then be nonzero, implying the existence of an edge $(V_j^{\alpha^{(j)}_{s}}, V_{j+1}^{\delta_{s}})$ in $G$ which does not lie in $C$. But by assumption $C$ contains all edges in $G$.

By Eq.~\ref{eqpt31}, $\mathrm{sign}(T) = (-1)^{N_e}$, where $N_e$ is the number of e-cycles in the subgraph associated with $T$. Since this subgraph is precisely $C$, $N_e = 1$ and $\mathrm{sign}(T) = -1$. Thus $A[\alpha^{(0)}] < 0$ and $A$ fails to be a $P_0$-matrix. \qquad \qed
\end{pf}

Corollary~\ref{corconv} illustrates an application of Lemma~\ref{lemconv}:

\begin{cor}
\label{corconv}
Suppose a graph $G_{A^{(0)}A^{(1)}\cdots A^{(k-1)}}$ contains an e-cycle $C$. Then there are matrices in $\mathcal{X} \equiv \mathcal{Q}(A^{(0)})\mathcal{Q}(A^{(1)})\cdots \mathcal{Q}(A^{(k-1)})$ which are not $P_0$-matrices. 
\end{cor}

\begin{pf}
Each edge in $C$ corresponds to an entry in one of the matrices $A^{(j)}$. For each $j$ define $\tilde A^{(j)}$ to be the matrix $A^{(j)}$ with all entries {\em not} corresponding to edges in $C$ set to be zero. Then the matrix factorisation $\tilde A = \tilde A^{(0)}\tilde A^{(1)}\cdots \tilde A^{(k-1)}$ gives rise to a graph which consists solely of the e-cycle $C$, and hence, by Lemma~\ref{lemconv}, $\tilde A$ fails to be a $P_0$-matrix. But $\tilde A \in \mathrm{cl}(\mathcal{X})$ (that is the closure of $\mathcal{X}$), and since the set of $P_0$-matrices is closed, there are matrices in $\mathcal{X}$ which fail to be $P_0$. \qquad \qed 
\end{pf}

\begin{pot4}
This is immediate: Theorem~\ref{corconv1} is simply a combination of Corollaries~\ref{cormain}~and~\ref{corconv}. \qquad \qed
\end{pot4}

\section{Notes and conclusions}

Although the special case $k=1$ was treated for completeness, the result can easily be inferred from previous work. For $k=1$, Theorem~\ref{corconv1} states that given a square matrix $A$, all matrices in $\mathcal{Q}(A)$ are $P_0$-matrices if and only if the (unique) signed $(1,\{1\})$-BC digraph $G_A$ associated with $A$ is e-cycle-free. However $G_A$ is closely related to an object often called the interaction graph or I-graph in the literature. In fact the I-graph associated with $A$ is just $G_{A^T}$. Results in \cite{soule,banajiJMAA} showed that $G_{-A^T}$ lacks positive cycles if and only if all matrices in $\mathcal{Q}(A)$ are $P_0$-matrices. Trivially, $G_{-A^T}$ lacks positive cycles if and only if $G_{-A}$ lacks positive cycles. The definitions imply that e-cycles (resp. o-cycles) in $G_A$ are in one-to-one correspondence with positive cycles (resp. negative cycles) in $G_{-A}$. So $G_{-A}$ lacks positive cycles if and only if $G_A$ is e-cycle-free. Together these observations imply that matrices $\mathcal{Q}(A)$ are all $P_0$-matrices if and only if $G_A$ is e-cycle-free.

The case $k=2$ has also effectively been treated previously in \cite{banajicraciun2,banajiJMAA}, where the associated graphs were termed ``DSR graphs''. The main differences between the definition of a DSR graph in \cite{banajiJMAA}, and a signed $(2,\{1\})$-BC digraph here, are (i) directions on all edges are reversed, (ii) here, edge-labels have been ignored, while some computations on DSR graphs in \cite{banajicraciun2,banajiJMAA} involved edge-labels, and (iii) in the construction of the DSR graph a pair of identically signed edges $(\tilde v, v)$ and $(v, \tilde v)$ are replaced with a single undirected edge, with a view to removing o-cycles of length $2$ from the graph, and thus simplifying computation. This process neither creates nor destroys e-cycles, and so does not change the key fact that an absence of e-cycles implies that associated matrices are $P_0$-matrices. 

The treatment in \cite{banajicraciun2,banajiJMAA} also suggests that extensions obtaining sharper results by introducing edge-labelling and more complex computations on the graphs are possible. The most useful forms that such extensions might take depend on the applications in question. These directions will be treated in future work.

\section*{References}

\bibliographystyle{unsrt}

\end{document}